\documentclass[12pt]{article}

\usepackage{graphicx}
\usepackage{amssymb}
\usepackage{amsmath}
\usepackage{amscd}
\usepackage{latexsym}
\usepackage{url}
\newcommand{\be}{\begin{equation}}
\newcommand{\ee}{\end{equation}}

\newcommand{\dlt}{\delta}
\newcommand{\bt}{\beta}

\newcommand{\vp}{\varphi}

\newcommand{\al}{\alpha}
\newcommand{\ra}{\rightarrow}
\newcommand{\sgm}{\sigma}

\newcommand{\gm}{\gamma}

\newcommand{\lgl}{\langle}
\newcommand{\rgl}{\rangle}
\newcommand{\prt}{\partial}

\begin{document}

\begin{center}

{\Large {\bf A Resolution of St. Petersburg Paradox} \\ [5mm]

V.I. Yukalov } \\ [3mm]

{\it
Bogolubov Laboratory of Theoretical Physics, JINR, Dubna 141980, Russia} \\
E-mail: yukalov@theor.jinr.ru

\end{center}

\vskip 1cm

\begin{abstract}

The St. Petersburg paradox is the oldest paradox in decision theory and has played
a pivotal role in the introduction of increasing concave utility functions embodying
risk aversion and decreasing marginal utility of gains. All attempts to resolve it
have considered some variants of the original set-up, but the original paradox has
remained unresolved, while the proposed variants have introduced new complications
and problems. Here a rigorous mathematical resolution of the St. Petersburg paradox
is suggested based on a probabilistic approach to decision theory.

\end{abstract}

\vskip 2cm

{\parindent=0pt
\noindent
\vspace{0.5cm}


\vskip 5mm
Declarations of interest: none }

\newpage

\section{Introduction}

The famous St. Petersburg paradox is, probably, the oldest paradox in decision theory,
and can be said to have promoted the birth of modern decision theory itself. It was
invented by Nicolas Bernoulli in 1713 and discussed in his private letters. The paradox
was formally stated by his cousin Daniel Bernoulli (1738), who worked in St. Petersburg.
The paradox and the suggested solution were published in the Proceedings of the Imperial
Academy of Sciences of St. Petersburg (Bernoulli (1738)).

In economics, the St. Petersburg paradox has played a particularly important role in
pointing out situations in which supposedly rational decisions based on expected gains
or even expected increasing utilities are not endorsed by real rational human decision
makers. The St. Petersburg paradox has opened a flood of attempts to solve it, which
turn out all to modify it in one way or another. The most important change involves the
introduction of a concave utility function which, in the words of Cramer (1728), captures
the concept that ``men of good sense estimate money in proportion to the usage that they
may make of it'' (and not necessarily in proportion to its quantity). Motivated by the
St. Petersburg paradox, the introduction of concave utility functions, which embody risk
aversion and decreasing marginal utility of gains, remains the central pillar of modern
economic theory.

However, it turns out that the solution in terms of utility functions of Cramer (1728)
and Bernoulli (1738) is not completely satisfactory since, as is stressed below, a
slight change in the formulation makes the paradox reappear. This problem and the many
proposed solutions have been revisited recently by Clark (2002) and by Martin (2008).
Thus, Clark states "This seems to be one of those paradoxes we have to swallow". Or,
as Martin concludes "The St. Petersburg result is strange... The appropriate reaction
might just be to try to accept the strange result". As the suggested solutions all
modify the initial problem, it is important to emphasize that the original paradox has
remained unresolved.

Here, we propose a resolution of the St. Petersburg paradox based on the introduction 
of stochasticity of the decision making problem to address the paradox. Stochasticity 
is argued to be the natural intrinsic feature of human behavior.

To be precise, first of all, in Section 2, we identify the origin of the St. Petersburg
paradox. In Section 3, we briefly summarize the most known attempts to avoid it by
changing the original formulation and explain why such attempts do not resolve the
paradox in its general form. In Section 4, we briefly mention the models used in
stochastic decision making and discuss the justification of stochastic approach. Then,
in Section 5, we formulate the probabilistic approach to decision theory to be used 
and in Section 6 show how the paradox becomes resolved in the frame of this approach. 
In the probabilistic approach there exists a belief (disbelief) parameter whose 
estimation for the case of the St. Petersburg paradox is considered in Section 7. 
The role of the trial prior probability is discussed in Section 8. The existence of 
repeated Bernoulli games is mentioned in Section 9, and in Section 10, the so-called 
inverse St. Petersburg paradox is investigated. Section 11 concludes.

\section{Origin of St. Petersburg Paradox}

First, we have to correctly formulate the problem in order to understand where the
paradox is located. We emphasize that the game, which the paradox stems from, is a
correctly posed mathematical problem. To resolve a problem does not mean to change
it, but this means to find a week point in its attempted solution and to change the 
solution. The essence of the St. Petersburg paradox is not in the formulation of the 
game, but in the definition of what has to be accepted as an optimal lottery.

In the St. Petersburg paradox, one meets the possibility of infinite number of
actions, such as coin tosses. Dealing with infinity in applied sciences one keeps in 
mind the notion of potential infinity that is assumed to be defined as a limit of a 
sequence. Potential infinity is treated not as an actual number, but as a trend dealing 
with asymptotically large numbers, which are such large numbers that, though being 
finite, are greater than any number given apriori 
(Poincar\'e (1902), Moore (1990), Rucker (1995)).

Below, following Samuelson (1960, 1977), we accept this understanding of an infinite 
Bernoulli game as a limit of a sequence of finite games and consider the sequence of 
lotteries. Only then it is possible to mathematically correctly formalize how 
a person makes a decision to play a toss game and how much the game would be worth.

\subsection{Original Bernoulli game}

It is necessary to distinguish the Bernoulli game lying at the heart of the
St. Petersburg paradox and the paradox as such.

The game, as it has been originally formulated by Bernoulli (1738), is defined as
follows. Paying a fixed fee to enter the game, a player tosses a fair coin repeatedly
until tails first appear, ending the game. If tails appear for the first time at the
$n$-th toss, the player wins $x_n=2^n$ monetary units. Thus, the player wins $2$
monetary units if tails appear on the first toss, $2^2$ units, if on the second,
$2^3$ units, if on the third, and so on. The gain of $2^n$ units occurs therefore
with the probability $1/2^n$. The cost of the ticket to enter the game is proportional
to the number of allowed tosses. The main question is: How much should a player pay
for the ticket to enter the game? Or in other words: How many tosses should a player
buy?

\vskip 2mm
{\bf Remark}. The nature of the monetary units is not important. It could be Dollars,
Euros, Francs, or any other units. Therefore, for the sake of brevity, we shall omit in
what follows the use of the term "monetary units", when this cannot lead to confusion.

\vskip 2mm

In mathematical parlance, the above game reads as follows (Samuelson (1960, 1977)).
There is a set of lotteries $\{L_n:\; n = 1,2,\ldots,\infty \}$ labeled according to
the number of allowed tosses, so that the first lottery is
$$
L_1 =
\left\{ 2 , \; \frac{1}{2} ~ \big| ~ 0 , \; \frac{1}{2} \right\} \; ,
$$
the second lottery is
$$
L_2 =
\left\{ 2 , \; \frac{1}{2} ~ \big| ~ 2^2 , \; \frac{1}{2^2} ~ \big|
~ 0 , \; \frac{1}{2^2}   \right\} \;  ,
$$
and so on, with the $n$-th lottery being
\be
\label{1}
L_n =
\left\{ 2 , \; \frac{1}{2} ~ \big| ~ 2^2 , \; \frac{1}{2^2} ~ \big| \ldots
\big| ~ 2^n , \frac{1}{2^n} ~ \big| ~ 0 , \; \frac{1}{2^n} \right\} \; ,
\ee
while, as has been stressed above, the lottery $L_\infty$ implies the limit of the 
series $\{L_n\}$ for $n \ra \infty$. 

The probability distribution of each lottery satisfies the normalization condition
$$
 \sum_{m=1}^n \frac{1}{2^m} + \frac{1}{2^n} = 1 \;  .
$$
Then the principal question "for how many tosses the player should buy a ticket" is
equivalent to the question "which lottery can be treated as optimal in the sense of
being preferred by the players?"

The above game is a well formulated mathematical problem and does not contain any
paradox in its formulation.

\subsection{Original St. Petersburg paradox}

The paradox appears in the attempt to suggest a normative mathematically justified way
of choosing an optimal lottery from the given set of lotteries in the Bernoulli game,
so that the prescribed choice would be in agreement with the choice of real players.

In the original setup, one compares the lottery expected utilities. The lottery enjoying
the largest expected utility is nominated {\it optimal}. One considers a linear utility
function $u(x) = x$ of payoffs $x$. In general, one could use an affine utility function
as $u(x) = const + c x$. It is clear, however, that without loss of generality, the
constant term can be omitted and the coefficient $c$ can be set to one by choosing the
appropriate monetary units. The St. Petersburg paradox does not depend on these minor
details. So, for the simplicity of notation, it is sufficient to use the linear function.

The lottery expected utility is commonly considered as the lottery value or the price
of the ticket that is worth paying for entering the lottery. It is exactly this meaning
that is kept in mind, when one talks of the lottery value in the St. Petersburg paradox.
This value directly depends on the prize that one can win by playing the game represented
by the lottery. Therefore, one usually considers as equivalent the terms lottery utility,
lottery price, and lottery payoff.

Thus, using the rules of expected utility theory (von Neumann and Morgenstern (1944),
Savage (1954)), one comes to the utility of lottery (\ref{1})
\be
\label{2}
 U_n = \sum_{m=1}^n u\left(2^m\right) \; \frac{1}{2^m} = n \;  .
\ee
As is evident, the largest expected utility corresponds to the lottery with the
infinite number of tosses:
\be
\label{3}
U_{max} = U_\infty \qquad ( n \ra \infty) \; .
\ee
Since the largest expected utility, with infinite payoff, assuming the infinite
number of tosses, by definition is optimal, one is led to conclude that any rational
player should decide to always enter the game for any arbitrarily large price. However,
in reality, real human players would not play for more than a few monetary units
(Bernoulli (1738)), hence the paradox.

\subsection{General formulation of Bernoulli game}

It is clear that the considered problem can be generalized to a set $\{L_n\}$ of the
lotteries
\be
\label{4}
 L_n = \{ x_1,p_1~ | ~ x_2, p_2 ~ | ~ \ldots ~ | x_n, p_n ~ | ~ 0 , p_n \} \; ,
\ee
with payoffs $x_m$ and their associated probabilities $p_m$. These probabilities can
be either objective (von Neumann and Morgenstern (1944)) or subjective (Savage (1954)),
as this has no impact on the nature of the problem. But, of course, the probability
measure over the outcomes is required to be normalized,
\be
\label{5}
 \sum_{m=1}^n p_m + p_n = 1 \;  .
\ee
The question is: Which lottery from the set of lotteries
$\{L_n:\; n = 1,2,\ldots,\infty\}$ is optimal?

\subsection{General formulation of St. Petersburg paradox}

According to utility theory, the lottery with the largest expected utility
\be
\label{6}
U_{max} = \max_n U_n \; , \qquad U_n = \sum_{m=1}^n u(x_m) p_m \; ,
\ee
is accepted as optimal. However, it is easy to show (Seidl (2013)) the following.

\vskip 2mm
{\bf Theorem 1}. For each nonincreasing probability distribution $p_m$, there exist
nondecreasing utility functions $u(x_m)$ such that the expected utility $U_n$ diverges
as $n\ra\infty$.

\vskip 2mm
{\it Proof}. For each probability distribution $p_m$, it is possible to define such
utility functions $u(x_m)$ that for all $m > m^*$ one has the inequality
$$
\frac{u(x_{m+1})p_{m+1}}{u(x_m)p_m} \geq 1 \qquad ( m \geq m^* ) \; .
$$
Then by the d'Alembert ratio test the sequence $\{U_n\}$ diverges as $n\ra\infty$. $\square$

\vskip 2mm
As an explicit example, let us take the linear utility function and the standard for
the St. Petersburg paradox probability distribution $p_m = 1/2^m$, then
\be
\label{7}
U_n = \frac{x(2^n-x^n)}{2^n(2-x)} \qquad
\left( u(x_m) = x^m , ~ p_m = \frac{1}{2^m} \right) \; .
\ee
In particular, for $x = 2$ we have
\be
\label{8}
U_n = n \qquad ( x = 2 ) \;  .
\ee
Therefore $U_n$ diverges, as $n \ra \infty$, for $x \geq 2$ and converges for
$x < 2$.

If the optimal lottery, with the largest expected utility, corresponds to $n\ra\infty$
and has the infinite value, then any rational player should feel it profitable to spend
all available money to buy a ticket allowing for the maximal possible number of tosses.
But numerous empirical data drastically contradict this conclusion, since the majority
of real players prefer the lotteries with quite a modest number of tosses (Bottom et al.
(1989), Rivero et al. (1990), Vivian (2004), Hayden and Platt (2009), Cox et al. (2009),
Neugebauer (2010), Cox et al. (2018), Nobandegani and Shultz (2020)). This contradiction
is the essence of the paradox.

\section{Attempts of Paradox Solution}

The contradiction between the theoretical prediction of what should be the optimal
lottery and the understanding of real human beings, appearing in the St. Petersburg
paradox, looks very strange. Several attempts have been made to resolve the paradox.
However all these attempts have not addressed the original paradox, but suggested
changes to the original formulation by replacing it with some other problems. Such
replacements cannot be considered as providing solutions of the original paradox
(Clark (2002), Martin (2008)).

\subsection{Use concave utility functions}

The first suggestion was made by Cramer (1728) and Bernoulli (1738) who proposed to
take a utility function characterized by a decreasing marginal utility, such as the
square-root function (Cramer (1728)) or the logarithmic utility function $u(x)=\ln x$
(Bernoulli (1738)), leading to finite expected utilities. Taking, for instance, the
logarithmic utility function gives
$$
\sum_{m=1}^\infty \; \frac{m\ln 2}{2^m} =  2\ln 2 \; .
$$
This approach opened the road towards the classical decision theory under uncertainty,
in which risk aversion and decreasing marginal utility of gains are quantified by
increasing concave utility functions.

However, this so-called ``solution" does not solve the initial paradox, since the
utility function has been modified, in the sense that it still accepts that a supposedly
well-defined lottery with the linear utility function may deliver an infinite average
gain. Furthermore, even if we accept the change of utility function, it is straightforward
to recover a paradox of the same type by changing the prize structure. Suppose, for
instance, that, instead of paying $x_m = 2^m$ for a run of $m$ heads followed by tails,
the prize is set to $x_m=\exp (2^m)$. Then, the expected utility, with the logarithmic
utility function, becomes infinite again. As suggested by Menger (1967), it is possible
to generalize this result and state that, for any utility function, there is a prize
structure such that the expected utility diverges.

\subsection{Assume the boundness of utility functions}

Menger (1967) argued that the only way of resolving the paradox is to assume that
there is an upper limit to utility. The same assumption was accepted by Hardin (1982),
Jeffrey (1983), and Gustason (1994). In that case, the expected utility remains finite
whatever the prize structure.

The finiteness of the utility function finds a justification in the neuroeconomics of
sensation. There are arguments that, because of the finite rates of neuronal firing
of action potentials, subjective values should be non-increasing beyond some point 
(Rayo and Becker (2007), Glimcher and Tymula (2018)). A similar saturation of sensation 
occurs for major disasters for which human moral intuitions are not well-calibrated to 
the extreme scales of deaths. This is called moral numbing (Slovic et al. (2020)).

Martin (2008) argued that limiting the utility function by some value is a quite 
subjective procedure. Such limits, if any, would be very different for different people 
and for different circumstances. They certainly cannot be applied to all people without 
exception with respect to money. There exist a number of examples when for people there 
is no limit of money utility, and the more they have, the better (Martin (2008)).

\subsection{Change the expected utility definition}

Generally, it is always possible to change the definition of the expected utility in
order to obtain a finite value for the newly defined quantity. For instance, one can
include some risk-aversion factors in the definition (Weirich (1984)).

There are however several objections to this approach. In fact, some people are not
risk averse in certain situations and can even enjoy taking risk. It is also possible
to adjust the prize structure in order to compensate the risk aversion and make again
infinite the expected utility, as discussed by Martin (2008). Using cumulative prospect
theory (Kahneman and Tversky (1979)), instead of expected utility theory, also does not
remove the paradox, unless one artificially selects some special probability distributions
and parameters (Rieger and Wang (2006)), which transforms the problem into a form having
little relation to the original Bernoulli setup. Blavatskyy (2005) showed that the 
cumulative prospect theory cannot explain the paradox, since it requires to overweight 
small probabilities which more than offsets the concavity of wealth. In any case, changing 
the rules of the game is not a solution to the original problem. As Pfiffermann (2011) 
concludes ``Cumulative prospect theory does not explain the St. Petersburg paradox".

\subsection{Set a finite limit to prizes}

Since no individual and no casino possess infinite resources, the total prize should
be limited by a finite value (Samuelson (1960), Gustason (1994)). This means that the
summation in $U_\infty$ should be cut at a finite $N$ by assuming that, either the
utility function or the probabilities take zero values for all $n > N$.

Mathematically, this corresponds to changing the original Bernoulli game, by either a
different definition of the utility function or of the payoffs, or of the probabilities.
Therefore, it cannot be considered to be a solution to the initial problem.

\subsection{Set a finite number of tosses}

Many argued (see, e.g., Nathan (1984) and references therein) that, in reality, nobody
can toss a coin an infinite number of times, so that the summation in $U_\infty$ must
be limited by a finite $N$.

Mathematically, this is just a particular case of the same modification of the game
as described in the previous variant. Therefore, the original problem remains open.
As Martin (2008) says, "refusing to think about the problem is not solving it".

\subsection{Replace the single game by a set of repeated games}

Sometimes one replaces the single Bernoulli game described above by a set of repeated
games played $N$ times, as is described by Vivian (Vivian (2013)). The set of $N$ 
repeated games can be treated as a compound lottery characterized by an expected monetary
value $U_N$ defined as the average per game monetary value. As has been found 
(Vivian (2013)), this expected value is $U_N = 1+\log_2 N$. Although $U_N$ increases
to infinity with $N \ra \infty$, but this increase is much slower than the rise of 
$U_n = n$ with $n \ra \infty$ in the original Bernoulli game. Experimental data confirm 
this behaviour (Klyve and Lauren (2011), Vivian (2013)). 

However, replacing one game by another is not the solution of the problem. Moreover, from 
the mathematical point of view the same problem remains. In the original Bernoulli game,
one is given a set of alternatives valued as $U_n = 1,2,\dots,\infty$, and the question is 
posed of deciding which of the alternatives is preferable. If the criterion of preference 
is the maximal $U_n$, then the choice is $U_\infty$. But the willingness to pay is finite 
and not large.

In the same way, if one is given a collection of expected values $U_N = 1 + \log_2 N$, 
with $N = 1,2,\ldots,\infty$, and the choice is prescribed by the largest $U_N$, then this
choice is again $U_\infty$, while the willingness to pay is finite. As earlier, the problem 
is not in the formulation of the game, but in the contradiction between the largest 
expected value and the willingness to pay. A detailed consideration is given in Sec. 9.

\subsection{Conclusion on the attempted ``solutions"}

By changing the rules of the Bernoulli game, it is certainly possible to adjust
them so that the expected utility or an equivalent decision functional be finite.
All previously suggested ``solutions" are exactly of this type. One solves some
variants that are not the original problem. In that sense, many versions could be
invented to replace the original St. Petersburg paradox (Samuelson (1977)). However,
``imposing restrictions on theory to rule out St. Petersburg bath water would throw out
some babies as well" (Martin (2008)). The commonly accepted conclusion is that ``from
the mathematical and logical point of view, the  St. Petersburg paradox is impeccable"
(Resnik (1987)). The suggested ``solutions" are not mathematical (Seidle (2013)). That
is, the original St. Petersburg paradox remains unresolved.

\section{Brief Account on Stochastic Decision Making}

The resolution of the St. Petersburg paradox, advanced in the present paper, is based 
on a probabilistic approach in decision theory. Therefore it is necessary to mention,
at least briefly, the main models of this approach and to emphasize the basic points 
justifying it.  

Generally speaking, in decision theory there are several approaches involving
probabilistic choice. The simplest model is a tremble model of Harless and Camerer
(1994). According to this model, individuals have a unique preference relation on
the set of lotteries, but they do not choose following their preferences all the time.
Instead, because of a random involutionary tremble occurring with some probability $p$
individuals choose an alternative which is not their preferred option. Individuals act
in accordance with their preferences with probability $1-p$. However, Carbone (1997)
and  Loomes et al. (2002) find that this tremble model fails to explain experimental
data.

The other is the Fechner (1860) model of random errors where individuals possess a
unique preference relation on the set of lotteries, but they reveal their preferences
with a random error as a result of carelessness, slips, insufficient motivation etc.
In the Fechner model, random error is usually drawn from a normal distribution with
zero mean and constant standard deviation (Hey and Orme (1994)) or from a truncated
normal distribution (Blavatskyy (2006)). Under the Fechner model, errors are more likely 
when the utilities of the options are close and less likely when they are less close. 
That is not the same as in the tremble model, where every choice over lottery pairs has 
the same error rate.

The third model of stochastic choice is a random preference model, where it is assumed
that individuals are endowed with several preference relations over the set of risky
lotteries and a probability measure over these preference relations. When faced with
a decision problem, individuals first draw a preference relation and then choose an
alternative which they prefer according to the selected preference relation (Loomes and
Sugden (1995)).

Finally, in the random choice model one assumes that individuals possess a probability
measure that captures the likelihood of a lottery being chosen over other lotteries.
Then there is no need for a mediator between a deterministic preference relation and an 
empirical stochastic choice pattern. The lotteries $L_n$ and, respectively, the related 
expected utilities are considered as random quantities characterized by a probability 
distribution $p(L_n)$ over the lotteries (Luce (1959, 1977), Gul et al. (2014)). One 
should not confuse the probabilities over the set of lotteries with the probabilities 
of payoffs inside each lottery.   

The probabilistic approach we employ in the present paper is based on the ideas of the
random choice model, however with an important generalization. The probability distribution
over the set of lotteries is not just postulated but it is defined by the principle of 
minimal information taken in the form of the Kullback-Leibler information functional 
(Kullback and Leibler (1951), Kullback (1959)), where the Luce rule plays the role of 
a trial prior. The method of deriving probability distributions from the minimization of 
an information functional, or from the minimization of a conditional relative entropy 
under given constraints, is well known and widely used. The most mathematically grounded
approach is based on the Shore-Jonson theorem (Shore and Jonson (1980)) stating that, 
given a prior density and new constraints, there is only one posterior density satisfying
these constraints and the conditions of uniqueness, coordinate invariance, and system
independence, such that this unique posterior can be obtained by minimizing the
Kullback–Leibler information functional (Kullback and Leibler (1951), Kullback (1959)).

Numerous empirical studies demonstrate that practically all decision makers do not 
follow the absolutely deterministic choice between the given alternatives, but the 
decisions vary for different subjects and even for the same subject at different 
moments of time (Luce (1958, 1959, 1977), French and Insua (2000), Gul et al. (2014)).

The most important is that stochasticity in decisions is unavoidable even for a single 
decision maker in the process of each decision. This is due to the fact that, as has 
been discovered in neurophysiological studies, the decision process in the brain is 
accompanied by noise leading to random choice, with the randomness caused by generic 
variability and local instability of neural networks (Werner and Mountcastle (1963), 
Arieli et al. (1986), Gold and Shadlen (2001), Glimcher (2005), Schumacher et al. (2011), 
Shadlen and Shohamy (2016), Webb, R. (2019), Kurtz-David et al. (2019)). Neurophysiological 
studies provide evidence that the process of decision making in regions of the human brain 
that are known to be responsible for value estimation is highly variable. This variability 
causes choice inconsistency at repeated decisions. The variability arises not merely from 
limitations in the data available to the decision maker, but the choice itself is stochastic 
due to fundamental stochasticity in the physiological process of making decisions 
(Werner and Mountcastle (1963), Arieli et al. (1986), Gold and Shadlen (2001), 
Glimcher (2005), Schumacher et al. (2011), Shadlen and Shohamy (2016), Webb, R. (2019), 
Kurtz-David et al. (2019)). 

In the recent review, Woodford (2020) summarizes the modern point of view that considers 
randomness as an internal feature of functioning of the human brain, where decisions are 
formed on the basis of noisy internal representation. The cognitive process through which 
judgments are generated involves two stages: encoding of the stimulus, that is the process 
through which the internal representation is produced, and decoding of the internal 
representation about the stimulus, both these stages are random, being accompanied by noise. 
The existence of random noise is caused by the way that neurons fire in response to the 
signals they receive. 

In this way, psychophysical and neurophysiological studies evidence that cognitive 
imprecision is due to the inevitable internal noise. This noise in the nervous 
system corrupts the evidence that must subsequently be decoded to produce a judgment.
Stochasticity is the unavoidable feature of the human brain functioning. As a result, 
the choice in decision making is not deterministic, being based on the comparison of 
utilities, but it is rather stochastic and based on the comparison of probabilities. 
It is important to stress that stochasticity is an internal property of a subject and 
it is always involved in the process of decision making at both the stages of encoding 
external signals and decoding them for formulating a decision. The internal stochastic 
fluctuations in the brain, which can be connected with the characteristics of external 
signals, are revealed in decision making as hesitations.

In the case of a subject deciding to take part in the Bernoulli gamble, hesitations can 
arise with respect to the rules of the gamble. First of all, a subject can have doubt
whether he/she has correctly understood the gamble description. The possibility of gaining
an enormous, or even infinite, payoff should certainly be suspicious for any reasonable
person. A normal human, when assuming the possibility of an infinite gain would for sure
feel that there is a kind of cheating here. Moreover, average humans have some limit
of imagination with respect to the amount of money, above which numbers sound as purely 
fantastic. Such suspicions should inevitably induce strong disbelief with respect to the 
gamble. In turn, this disbelief would amplify the naturally existing stochasticity in the 
internal encoding--decoding process in the brain. 

This is why in order to understand the origin of the St. Petersburg paradox, it is 
necessary to resort to the probabilistic approach in decision making. The Bernoulli game 
as such is absolutely correct and does contain an infinite lottery with infinite expected
utility. However, for real humans the optimality of a lottery is not prescribed by its 
largest expected utility, but it is defined by the largest probability of the lottery. 
So there is no need to change anything in the Bernoulli game that is mathematically 
precise. The paradox is not in the game, but in the application to the game of the 
deterministic preference criterion. Therefore the paradox can be eliminated only by
using the probabilistic preference criterion.

\section{Formulation of the Approach}

Suppose that a set of alternatives is presented as a set of lotteries 
$\{L_n: n =1,2,\ldots\}$. For the sequence $\{L_n: n =1,2,\ldots\}$ of these lotteries, 
there corresponds the sequence $\{U_n \equiv U(L_n): n =1,2,\ldots\}$ of their expected 
utilities. The notion of the infinite lottery is defined as the limit of the sequence 
$\{L_n\}$ for $n$ tending to infinity, $n \ra \infty$. According to the random choice 
approach, the lotteries are treated as random variables, with a probability measure 
$\{p(L_n): n=1,2,\ldots\}$ over the set of these lotteries. In other words, the choice 
is not deterministic, being based on the comparison of utilities, but it is rather 
stochastic and based on the comparison of probabilities. Below we formulate these 
properties explicitly. The use of stochastic preferences is the basic point of the 
approach based on the following definitions.   

\vskip 2mm

{\bf Definition 1}. The lottery $L_i$ is stochastically preferred to $L_j$ if and only
if
\be
\label{9}
p(L_i) > p(L_j) \qquad ( L_i \succ L_j ) \; .
\ee

\vskip 2mm
{\bf Definition 2}. The lotteries $L_i$ and $L_j$ are stochastically indifferent if and
only if
\be
\label{10}
p(L_i) = p(L_j) \qquad ( L_i \sim L_j ) \;   .
\ee

\vskip 2mm
{\bf Definition 3}. The lottery $L_{n_{opt}}$ is called stochastically optimal if and
only if
\be
\label{11}
 p(L_{n_{opt}}) = \sup_n p(L_n) \; .
\ee

\vskip 2mm

A stochastically optimal lottery corresponds to the maximal lottery probability. This 
implies that, according to the frequentist interpretation of probability, the majority 
of subjects choose exactly the stochastically optimal lottery where the probability 
possesses a maximum.  

\vskip 2mm

It is possible to study general properties of the probability measure and the process
of stochastic optimization (Rubinstein and Kroese (2004), Gul et al. (2014)). However
we concentrate our attention on the practical way of constructing an explicit expression
for the probability, since its explicit formula is compulsory for the resolution of the
St. Petersburg paradox.

The probability measure over the set of lotteries has to satisfy the necessary 
conditions (Gul et al. (2014), Yukalov and Sornette (2014, 2018)) characterizing the 
sought probability. The first condition is the standard normalization for the probability
\be
\label{12}
\sum_{n=1}^\infty \; p(L_n) = 1 \; , \qquad 0 \leq p(L_n) \leq 1 \;  .
\ee
The second condition is the existence of a global mean 
\be
\label{13}
\sum_{n=1}^\infty \; p(L_{n}) U_n = U  \; , \qquad (|U| < \infty) \; .
\ee
 
One should not confuse the probability $p_n$ of payoffs $x_n$ in a lottery with the 
probability of a lottery $p(L_n)$ as a whole. Respectively, one should not confuse
the expected utility $U_n$ of a lottery $L_n$ with the global mean $U$ that is not
an expected lottery utility, but just a constraint. The explicit form of the probability
distribution $p(L_n)$ can be derived from the minimization of the Kullback-Leibler 
information functional (Kullback and Leibler (1951), Kullback (1959), 
Shore and Jonson (1980)). The Shore-Jonson theorem (Shore and Jonson (1980)) states that, 
given a prior density and additional constraints, there is only one posterior density 
satisfying these constraints and the conditions of uniqueness, coordinate invariance, and 
system independence, such that this unique posterior can be obtained by minimizing the 
Kullback–Leibler information functional (Kullback and Leibler (1951), Kullback (1959)). 

It is important to stress the following pivotal point. If one postulates the explicit form 
of the Kullback–Leibler information functional and defines the probability $p(L_n)$ as its 
minimizer, then this makes a sufficient condition for deriving the expression of $p(L_n)$, 
without requiring the finiteness of $U$ in constraint (\ref{13}). However the Shore-Jonson 
theorem (Shore and Jonson (1980)) is more general. In addition to the above sufficient 
condition, it states that, under given constraints and conditions of uniqueness, coordinate
invariance, and system independence, the posterior probability $p(L_n)$ is necessarily 
the minimizer of a functional enjoying the Kullback–Leibler form, provided the imposed 
constraints, such as (\ref{12}) and (\ref{13}) do not contain singularities, which here 
means a not divergent value of $U$. The Shore-Jonson theorem justifies why it is precisely 
the Kullback–Leibler functional that has to be accepted for the minimization procedure. 

In the present case, this information functional for the posterior probability distribution 
$p(L_n)$, under a prior distribution $p_0(L_n)$ and constraints (\ref{12}) and (\ref{13}), 
reads as
\be
\label{14}
I[\; p \;] = \sum_{n=1}^\infty \; p(L_{n}) \ln \; \frac{p(L_n)}{p_0(L_n)} +
\al \left[ \; 1 - \sum_{n=1}^\infty \; p(L_{n}) \; \right] +
\bt  \left[ \; U - \sum_{n=1}^\infty \; p(L_{n}) U_n \; \right] \;  ,
\ee
where $\alpha$ and $\beta$ are Lagrange multipliers. The Kullback-Leibler information
functional is a measure for the amount of information necessary to transform a prior
probability distribution $p_0(L_n)$ into the posterior distribution $p(L_n)$. The choice
of a prior distribution will be explained below.

The minimization of the Kullback–Leibler information functional is equivalent to the 
minimization of the Kullback–Leibler relative entropy, or divergence, under given 
constraints. This concept is well known and widely employed in statistical models of 
inference, different branches of applied statistics, information theory, computing, 
machine learning, pattern recognition, development of artificial intelligence, mechanics, 
neuroscience, theory of biological systems, and various kinds of optimization problems 
(Kullback and Leibler (1951), Kullback (1959), Burnham and Anderson (2002), MacKay (2003),
Nielsen (2005), Cover and Thomas (2006)). As is explained above, the Shore-Jonson theorem  
(Shore and Jonson (1980)) justifies why the posterior distribution has to be defined as 
the minimizer of the Kullback-Leibler functional. When minimizing this functional, one can
meet a finite index set $n =1,2,\ldots,N$ as well as the infinite index set 
$n = 1,2,\ldots,\infty$. The rigorous procedure of the minimization is done for a finite 
index set, with a finite $N$, and the transition to an infinite index set 
$n = 1,2,\ldots,\infty$ is accomplished after the minimization procedure by sending 
$N \ra \infty$. The details of the transition $N \ra \infty$ can be found in literature 
cited above and in (Dmitruk and Kuzkina (2005), Pinski (2015)). The minimization theorem 
for a finite $N$ is formulated below.  

\vskip 2mm

{\bf Theorem 2}. The probability distribution $p(L_n):\; n = 1,2,\ldots,N$, defined as the 
minimizer of the Kullback-Leibler information functional 
$$
 I_N[\;p\;] = \sum_{i=1}^N p(L_n) \;\ln \;\frac{p(L_n)}{p_0(L_n)} \; + \;
\al \left[\; 1 - \sum_{n=1}^N p(L_n) \; \right] \; + \; 
\bt \left[\; U - \sum_{n=1}^N p(L_n)  U_n\; \right] \; ,
$$
reads as
\be
\label{15}
 p(L_n) =
\frac{p_0(L_n)e^{\bt U_n} }{\sum_{n=1}^\infty p_0(L_n)e^{\bt U_n} } \; .
\ee

\vskip 2mm

{\it Proof}. The proof is based on the procedure of minimizing functionals in calculus 
of variations (Gelfand and Fomin (1963)). The necessary condition for the functional 
$I[p]$ to be an extremum is the zero first-order variation
$$
 \dlt I[\;p\;] = \sum_n \frac{\prt I[\;p\;]}{\prt p(L_n)} \; \dlt p(L_n) = 0 \;  .
$$
From here, using the Kullback-Leibler information functional, we have
$$
\frac{\prt I[\;p\;]}{\prt p(L_n)} = \ln \;\frac{p(L_n)}{p_0(L_n)} + 
1 - \al - \bt U_n = 0 \;   ,
$$
which, together with the normalization condition, yields expression (\ref{15}). 

The necessary condition for the extremum to be a minimum is the non-negativity of the 
second-order variation. This condition becomes also sufficient if the second-order 
variation is strictly positive,
$$
 \dlt^2 I[\;p\;] = \frac{1}{2} \sum_{mn} I_{mn} \; \dlt p(L_m) \dlt p(L_n) > 0 \;  ,
$$
where 
$$
 I_{mn} = \frac{\prt^2 I[\;p\;]}{\prt p(L_m)\; \prt p(L_n)} \;  .
$$
The second-order variation is positive provided that the Hessian matrix $[I_{mn}]$, 
composed of the elements $I_{mn}$, is positive, which requires that all principal minors 
of this Hessian matrix be positive. It is easy to see that the Hessian matrix is diagonal, 
since
$$
I_{mn} = \frac{\dlt_{mn}}{p(L_n)} \;   .
$$
Therefore for its positivity, it is necessary and sufficient that 
$$
I_{nn} = \frac{\prt^2 I[\;p\;]}{\prt p(L_n)^2} = \frac{1}{p(L_n)} > 0  \;  ,
$$
that is the probability $p(L_n)$ has to be positive. $\square$ 
  
\vskip 2mm

{\bf Remark}. According to the standard arguments (Burnham and Anderson (2002), MacKay (2003),
Nielsen (2005), Cover and Thomas (2006), Dmitruk and Kuzkina (2005), Pinski (2015)),
after the minimization procedure, the finite index set can be extended to infinite by
sending $N$ to $\infty$. 

\vskip 2mm

As a trial prior, it is possible to take the Luce (1958, 1959, 1977) form defined
as follows. If the considered alternatives are characterized by the related attributes
$a_n \geq 0$, then the trial probability of choosing an alternative $L_n$ is
\be
\label{16}
 p_0(L_n) = \frac{a_n}{\sum_{n=1}^\infty a_n} \;  .
\ee
Hence the posterior distribution (\ref{15}) becomes
\be
\label{17}
 p(L_n) = \frac{a_n e^{\bt U_n}}{\sum_{n=1}^\infty a_ne^{\bt U_n} } \;  .
\ee
The first and direct attribute of a lottery is given by its expected utility. Therefore,
for the case of semi-positive expected utilities, the latter can be accepted as the
lottery attributes, while for negative expected utilities, the lottery attributes can
be defined as inverse lottery absolute values,
\begin{eqnarray}
\label{18}
a_n = \left\{ \begin{array}{cl}
   U_n ,          ~ & ~ U_n \geq 0 \\
|\; U_n\;|^{-1} , ~ & ~ U_n < 0
\end{array}
\right. \; .
\end{eqnarray}
Thus we obtain the explicit expressions for the probability distribution.

\vskip 2mm
{\bf Definition 4}. The probability distribution over the set of lotteries with
semi-positive expected utilities has the form
\be
\label{19}
p(L_n) = \frac{U_n e^{\bt U_n}}{\sum_{n=1}^\infty U_ne^{\bt U_n} } \qquad
( U_n \geq 0 ) \; ,
\ee
and for the lotteries with negative expected utilities,
\be
\label{20}
p(L_n) =
\frac{|\;U_n\;|^{-1} e^{\bt U_n}}{\sum_{n=1}^\infty |\;U_n\;|^{-1}e^{\bt U_n} }
\qquad  ( U_n < 0 ) \;   .
\ee

\vskip 2mm

It is useful to note that these distributions differ from the logit form, which is
often postulated. As is seen, the lottery probability $p(L_n)$ depends on the expected 
utility $U_n$ of this lottery. According to the criterion of stochastic preference, the
optimal lottery corresponds to the maximal lottery probability. 

\vskip 2mm
  
The parameter $\beta$, from the mathematical point of view, is a Lagrange multiplier 
guaranteeing the normalization condition (\ref{13}). From the psychological side, the 
parameter $\beta$ plays the role of a {\it belief parameter} reflecting the belief of 
a subject in the fairness of the gamble and in the subject confidence with respect to 
his/her understanding of the overall rules and conditions of the gamble. Under neutral 
belief, when $\beta = 0$, the probability distribution returns to the Luce form (\ref{16}). 
In the case of strong belief, the situation reduces to the deterministic choice of 
the largest expected utility:
\begin{eqnarray}
\label{21}
p(L_n) = \left\{ \begin{array}{ll}
   1 , ~ & ~ U_n = U_{max} \\
   0 , ~ & ~ U_n \neq U_{max}
\end{array} \right. \qquad ( \bt \ra \infty)  \; ,
\end{eqnarray}
where
\be
\label{22}
  U_{max} \equiv \max_n U_n \; .
\ee

It also may happen strong disbelief, when $\beta \ra -\infty$, and it becomes more
profitable to choose the lottery with the minimal expected utility. Although in realistic
cases, the belief parameter is finite. It can be positive in the case of belief or
negative in the case of disbelief. 

When all expected utilities are finite, that is $|U_n|<\infty$, then the probability 
distribution satisfies all normalization conditions by construction, including the 
normalization condition (\ref{13}). The situation is more restrictive, when the maximal 
expected utility (\ref{22}) diverges, so that $U_{max} \ra \infty$, as it happens in the 
St. Petersburg paradox.

\vskip 2mm
{\bf Theorem 3}. If the maximal expected utility (\ref{22}) diverges, so that
$U_{max}\ra\infty$, then the normalization condition (\ref{13}) requires that $\bt$ be
negative.

\vskip 2mm
{\it Proof}. From the expression for the probability distribution (\ref{19}) or
(\ref{20}), it follows that
\begin{eqnarray}
\label{23}
p(L_{max}) = \left\{ \begin{array}{ll}
   1 , ~ & ~ \bt \geq 0 \\
   0 , ~ & ~ \bt < 0
\end{array} \right. \qquad ( U_{max} \ra \infty)  \; .
\end{eqnarray}
If $\bt$ is semi-positive and $U_{max}\ra\infty$, then the normalization condition
(\ref{13}) does not hold, since $U\ra U_{max}\ra\infty$. This condition (\ref{13}) is
satisfied only for $\beta<0$. $\square$

\vskip 2mm

{\bf Remark}. Similarly, it is easy to show that if the minimal expected utility tends 
to $-\infty$, then $\bt$ must be semi-positive. In the case of the Bernoulli game, the
parameter $\bt$ has to be negative, since the expected utilities are positive defined
and can reach infinite value. A negative parameter $\bt$ signifies disbelief that can 
be caused e.g. by the disbelief of decision makers in the fairness of the game or by
their disbelief in the principal possibility of gaining an infinite sum.

\section{A Resolution of St. Petersburg Paradox}

Now we are in a position to formulate a resolution of the St. Petersburg paradox.
Recall that to resolve the paradox, it is necessary to explain, why in the Bernoulli
game, despite a divergent maximal expected utility, people consider as optimal a lottery
with a finite expected utility, and respectively, with a finite number of tosses.

\vskip 2mm
{\bf Theorem 4}. Assume that in the Bernoulli game the maximal expected utility diverges.
Then, under the probability distribution (\ref{19}), the stochastically optimal lottery
has a finite expected utility corresponding to a finite number of tosses.

\vskip 2mm
{\it Proof}. In the Bernoulli game, the expected utility of a lottery with $n$ tosses
is positive, $U_n>0$. By assumption, the largest expected utility diverges,
$U_{max}\ra\infty$, hence the belief parameter, according to Theorem 3, has to be
negative, $\bt<0$, corresponding to disbelief. Therefore the probability distribution
over lotteries has the form
\be
\label{24}
 p(L_n) = \frac{U_n}{Z} \; e^{\bt U_n} \;  ,
\ee
with the normalization factor
$$
Z = \sum_{n=1}^\infty U_n e^{\bt U_n} \qquad ( \bt < 0 ) \;  .
$$
By Definition 3 in equation (\ref{11}), the stochastically optimal lottery
$L_{n_{opt}}$ is defined by the maximum of the probability $p(L_n)$, where
\be
\label{25}
 \frac{\prt p(L_n)}{\prt U_n} = 0 \qquad ( L_n = L_{n_{opt}} ) \; .
\ee
Then the stochastically optimal lottery corresponds to the expected utility
\be
\label{26}
U_{opt} = U(L_{n_{opt}}) = \frac{1}{|\;\bt\;|} \; .
\ee
This utility is finite for any negative value of the belief parameter, and there 
is no paradox. $\square$

\vskip 2mm

For illustration, let us consider the original Bernoulli game, where $U_n=n$. This
leads to the probability
\be
\label{27}
 p(L_n) = \frac{n}{Z} \; e^{\bt n} \qquad ( \bt < 0 ) \;  ,
\ee
with the normalization factor
$$
Z = \sum_{n=1}^\infty n e^{-|\;\bt\;|n} = \frac{1}{4\sinh^2(|\;\bt\;|/2)} \;  .
$$

In the optimal lottery, the number of tosses is $n_{opt}\approx 1/|\bt|$. Keeping in
mind that $n_{opt}$ is to be a whole number, it lays in the interval
\be
\label{28}
 \left[ \; \frac{1}{|\;\bt\;|} \; \right] ~ \leq ~ n_{opt} ~ \leq ~
 \left[ \; \frac{1}{|\;\bt\;|} \; \right] + 1 \;  ,
\ee
where the square brackets imply the entier (the greatest integer not exceeding the
specified number).

Empirical investigations of the Bernoulli game (Cox et al. (2009, 2018), Hayden and
Platt (2009), Neugebauer (2010), Nobandegani and Shultz (2020)) show that the maximal
fraction of players prefers one or two steps in the game. Keeping in mind the frequentist
understanding of probability implies that the maximal number of players considers as
stochastically optimal the game of one - two tosses, that is $n_{opt} \approx 1, 2$.
Therefore, according to condition (\ref{28}), we have $|\beta| \approx 1$. Thus, despite
that the maximal expected utility of the Bernoulli toss game, is infinite, the majority
of real human beings would not evaluate it so much, but, in contrast, would not buy a
ticket for more than a few monetary units.

\section{Behavioral Meaning of Belief Parameter}

Mathematically, the parameter $\beta$ appears as a Lagrange multiplier in the process 
of the Kullback-Leibler information-functional minimization. Its interpretation as of 
a parameter characterizing the belief (or disbelief) of a decision maker is based on 
its role. It is termed {\it belief} parameter when $\beta > 0$ and {\it disbelief} 
parameter if $\beta < 0$. As is shown in Sec. 5, in the case of strong belief, when 
$\beta \ra \infty$, we return to the deterministic variant of decision theory prescribing 
that $100 \%$ of decision makers should prefer the lottery with the largest expected 
utility. In the opposite case of strong disbelief, when $\beta \ra -\infty$, one chooses 
the lottery with the minimal expected utility. And $\beta = 0$ characterizes the neutral 
situation, when there is neither noticeable belief nor disbelief.

In the considered case of the St. Petersburg paradox, the behavioral meaning of the 
belief parameter is justified by the attitude of a decision maker to the game. The 
formulation of the Bernoulli gamble suggests a seeming possibility of gaining infinite 
(or at least enormous) amount of money. For usual people, such a possibility sounds as 
highly unrealistic, because of which they do not give much trust to this promise. That is, 
a normal person does not believe in the fairness of the gamble. This is why the belief 
parameter is negative for this case.

As is evident, the value of the belief parameter is closely connected with the particular 
game. In the Bernoulli game, from the behavioural point of view, there is disbelief in
the fantastically sounding aptitude of winning infinite amount of money, because of which
this parameter is negative. Being a Lagrange multiplier, the value of this parameter,
as is usual for optimization problems, can be found from the knowledge of some 
phenomenological data or from additional constraints. For example, the value of the 
disbelief can be connected with the level of uncertainty in the game. For a set of 
lotteries, the typical quantity characterizing the level of uncertainty is the 
standard deviation (Ghahramani (2000)). Considering a set of lotteries $L_n$ with 
the expected utilities $U_n$ and probabilities $p(L_n)$, the standard deviation is
\be
\label{B1}
 \sgm =\sqrt{ {\rm var}(U_n) }  \; ,
\ee
which is a square root of the variance 
\be
\label{B2}
{\rm var}(U_n) = \lgl \; U_n^2 \; \rgl -  \lgl \; U_n \; \rgl^2 \; ,
\ee
where
\be
\label{B3}
 \lgl \; U_n \; \rgl = U = \sum_n U_n \; p(L_n) \; , \qquad
\lgl \; U_n^2 \; \rgl =  \sum_n U_n^2 \; p(L_n) \; .
\ee
The larger the standard deviation, the larger the uncertainty in the choice among the 
alternatives, hence the larger the absolute value of the disbelief represented by the 
parameter $\beta$. This allows us to consider the disbelief parameter being proportional
to the standard deviation,
\be
\label{B5}
  \bt = - |\;  \bt \; | = - \sgm \qquad ( \bt < 0 ) \;  .
\ee

For the original formulation of the Bernoulli game, we have the expected utilities 
$U_n=n$, with the probabilities (\ref{27}). Therefore we have
$$
\lgl \; U_n \; \rgl = U = \sum_{n=1}^\infty n \; p(L_n) = 
\coth\left( \frac{|\bt|}{2} \right) \; , 
$$
\be
\label{B6}
\lgl \; U_n^2 \; \rgl =  \sum_{n=1}^\infty n^2 \; p(L_n) = 
\frac{1}{2} \; \left( 3 U^2 - 1\right) \; .  
\ee
Then the variance (\ref{B2}) becomes
\be
\label{B7}
 {\rm var}(U_n) = \frac{1}{2} \; \left( U^2 - 1\right) = 
\frac{1}{2\sinh^2(|\;\bt\;|/2)} \;  .
\ee
In this way, for the disbelief parameter (\ref{B5}) we obtain the equation
\be
\label{B8}
\sqrt{2} \; |\bt | \sinh \; \frac{|\bt|}{2} = 1 \;  ,
\ee
whose solution gives $|\bt|=1.157$. This estimate, according to the previous 
section, tells us that the optimal lottery corresponds to $n_{opt}$ of order one,
that is, to the game of just a few tosses. This is in perfect agreement with the 
empirical observations that the majority of subjects prefer the Bernoulli game 
of one - two steps (Cox et al. (2009, 2018), Hayden and Platt (2009), Neugebauer (2010), 
Nobandegani and Shultz (2020)).

\section{Role of Prior Probability}

The form of the lottery probability $p(L_n)$ that is the minimizer of the information 
functional, generally, depends of the prior probability $p_0(L_n)$. The latter was 
chosen above in the simple form satisfying the Luce rule. This simple form has the 
advantage of containing no additional parameters complicating the consideration. 
However, it is possible to put the question: What are the general conditions imposed 
on the prior probability, which allow for the solution of the St. Petersburg paradox?

It is evident that the prior probability distribution, by its meaning, cannot be 
represented by an arbitrary distribution, but has to satisfy reasonable assumptions
appropriately describing the expected features of the problem. Thus, it is reasonable 
to assume that the prior probability, in addition to being semi-positive, should 
depend on the lottery expected utility $U_n$ so that, the increasing expected utility 
$U_n$ would increase $p_0(L_n)$, which means that
\be
\label{30}
 \frac{\prt p_0(L_n)}{\prt U_n} ~ > ~ 0 \;  .
\ee

Without the loss of generality, the prior probability can be written as
\be
\label{31}
p_0(L_n) = \frac{\vp(U_n)}{\sum_n \vp(U_n)}  \; .
\ee
The function $\vp$, in view of condition (\ref{30}), has the property
\be
\label{32}
 \frac{\prt \vp(U_n)}{\prt U_n} ~ > ~ 0 \; .
\ee
Then the minimizer of the information functional takes the form
\be
\label{33}
p(L_n) = \frac{\vp(U_n)\; e^{\bt U_n}}{\sum_n \vp(U_n)\; e^{\bt U_n}} \; .
\ee
   
The sufficient conditions for the existence of a finite stochastically preferred 
lottery, are the conditions of the probability maximum  
\be
\label{34}
 \frac{\prt p(L_n)}{\prt U_n} = 0 \; , \qquad
  \frac{\prt^2 p(L_n)}{\prt U_n^2} ~ < ~ 0 \qquad 
( U_n = U_{n_{opt}} ) \; .
\ee
By using the notation
\be
\label{35}
x \equiv U_n \; , \qquad x_{opt} \equiv U_{n_{opt} } \; ,
\ee
these conditions can be rewritten as
\be
\label{36}
\vp'(x) + \bt \vp(x) = 0 \; , \qquad  \vp''(x) - \bt^2 \vp(x) < 0
\qquad ( x = x_{opt} ) \;   ,
\ee
where the primes imply the derivatives over $x$. 

In this way, any function $\varphi$ satisfying conditions (\ref{32}) and (\ref{36}) 
will lead to the same general conclusion of the existence of a finite stochastically 
preferred lottery, that is, to the resolution of the paradox. Concrete formulas 
defining the optimal lottery will, of course, be slightly different. For instance, 
a simple function, satisfying condition (\ref{32}), could be
\be
\label{37}
 \vp(x) = x^\al \qquad ( \al > 0 ) \; .
\ee
Then the first of conditions (\ref{36}) gives 
\be
\label{38}
  x_{opt} = \frac{\al}{|\bt| } \; ,
\ee
while the second of conditions (\ref{36}) is valid for any $\al>0$. 

Another example is a logarithmic function $\varphi$, such that
\be
\label{39}
 \vp(x) = \ln \; ( 1 + x ) \; ,
\ee
with the notation
\be
\label{40}
x \equiv \frac{U_n}{U_0} \; , \qquad  
x_{opt} \equiv \frac{U_{n_{opt}} }{U_0} \;  ,
\ee
where $U_0 > 0$ is a scaling constant. Then the optimal lottery is defined by the 
equation
\be
\label{41}
 ( 1 + x_{opt} ) \; \ln \; ( 1 + x_{opt} )= \frac{1}{|\bt|} \;  ,
\ee
the second of conditions (\ref{36}) being always valid. 

One more example can be based on the general logit model that is widely employed in 
economics and stochastic decision theory. In the general version of the logit model,
the trial probability of an alternative $L_n$ is assumed to have the form
\be
\label{L1}
p_0(L_n) = \frac{e^{V_n}}{\sum_n e^{V_n}} \;  ,
\ee
where $V_n$ is the expected value of $L_n$. The expression for the value $V_n$ is to be
postulated and might be very different from expected utility, e.g. as in the stochastic
approach based on prospect theory (Kahneman and Tversky (1979)). In the majority of cases, 
the expected value $V_n$ is defined as a transformation $V_n = V(U_n)$ of the expected 
utility $U_n$. There exists a large variety of such transformations, including highly 
nonlinear ones (Small (1987), Gerken (1991), Wen and Koppelman (2001), Williams (2016)). 
For the illustration of the idea, let us take a weakly nonlinear transformation 
representing the value $V_n$ in the form
\be
\label{L2}
V_n = V(U_n) = b_n U_n^{\gm_n} + c_n \;   ,
\ee
where $b_n$, $c_n$, and $\gamma_n$ are positive parameters. These parameters are usually 
fitted to a particular problem of interest. As in the case of the St. Petersburg paradox,
the expected utility is assumed to be in the interval $U_n \in [0, \infty)$ and the belief
parameter $\beta$ is negative, expressing disbelief. The probability distribution (\ref{15}) 
becomes
\be
\label{L3}
p(L_n) = \frac{\exp(V_n + \bt U_n)}{\sum_n \exp(V_n+\bt U_n)} \;  .
\ee
The stochastically preferable alternative corresponds to the maximum of this probability,
which gives the alternative with the expected utility 
\be
\label{L4}
 U_{n_{opt}} = \left( \frac{b_n\gm_n}{|\;\bt\;|} \right)^{1/(1-\gm_n)} \;  ,
\ee
where $0 < \gamma_n < 1$. 

In the general case of a transformation $V_n = V(U_n)$, the stochastically preferable
alternative is defined by the equations
\be
\label{L5}
 \bt + \frac{\prt V_n}{\prt U_n} = 0 \; , \qquad 
\frac{\prt^2 V_n}{\prt U_n^2} < 0 \; .
\ee
        
Thus, there is a wide class of prior probabilities, satisfying straightforward 
conditions, which lead to the existence of a finite stochastically preferred lottery, 
hence resolving the paradox. It is worth emphasizing that this does not impose any 
constraints on the Bernoulli game as such, which does contain the limiting lottery 
with infinite expected utility. However the latter is not stochastically preferred.

\section{Repeated Games}

In the above consideration, we were keeping in mind a set of single Bernoulli games, 
consisting of $n = 1,2,\ldots$ tosses corresponding to the lotteries $L_n$ with the 
expected utilities $U_n$. The aim was to find an optimal lottery $L_{n_{opt}}$ that 
would be preferable. It is necessary to distinguish the expected utility or its 
particular form called the expected value, or the expected monetary value, and the 
willingness to pay. The expected value $U_n$ can increase to infinity with rising $n$, 
but the willingness to pay, that is denoted as $U_{n_{opt}}$, is finite. The mathematical 
formulation of the Bernoulli game is absolutely accurate and contains no paradoxes. 
Under the paradox, one assumes the apparent difference between the maximal expected value 
$U_{max} = \max_n U_n$ and the willingness to pay $U_{n_{opt}}$. In the probabilistic 
approach above, it is explained that these notions are really different and do not need 
to coincide, thus removing the seeming paradox. 

The same consideration is applicable to other games. For example, it is possible 
to consider a set of compound Bernoulli games, of $k$ tosses each, consisting of 
subsets of Bernoulli games repeated $N = 2^k$ times (Tversky and Bar-Hillel (1983), 
Klyve and Lauren (2011), Vivian (2013)). The expected value of a lottery is accepted 
as the average per game expected monetary value. The run of $N$ games corresponds to 
a compound lottery $L_N$, with the expected monetary value for all $N$ repeated games 
$N(1+\log_2 N)$, as is explained in the paper by Vivian (2013). This gives the average 
per game expected value  
\be
\label{42} 
 U_N = 1 + \log_2 \; N  
\ee
in some monetary units. This value $U_N$ increases with $N$, although much slower 
than $U_n=n$ in the single Bernoulli game (Tversky and Bar-Hillel (1983), Klyve 
and Lauren (2011), Vivian (2013)). Nevertheless, the expected value $U_N$ tends to 
infinity with $N$.
  
Similarly to the single game, we have to remember that the expected value $U_N$ does 
not need to equal the willingness to pay $U_{N_{opt}}$. Then one can ask the question: 
which of the runs would be preferred by decision makers, that is, what would be the 
willingness to pay? In other words, how many runs $N$ would be considered by subjects as 
optimal?     

Following the approach of Sec. 5, the probability of a lottery $L_N$ of $N$ repeated 
games is given by the expression 
\be
\label{43}
 p(L_N) = \frac{U_N \; e^{\bt U_N} }{\sum_{N=0}^\infty U_N \; e^{\bt U_N} } \; ,
\ee 
with the value $U_N$ defined in (\ref{42}). Here the parameter $\beta<0$ is different 
from that in probability (\ref{24}). The stochastically preferred $L_N$ is that which 
is prescribed by the maximum of probability (\ref{43}), corresponding to the subset
consisting of $N_{opt}$ repeated games. Therefore the willingness to pay is
\be
\label{44}
 U_{N_{opt}} = \frac{1}{|\;\bt\;|} \;  .
\ee
This defines the optimal number of repeated games by the equation
\be
\label{45}
 N_{opt} = 2^{1/|\bt|-1} \; .
\ee

Overall, the situation is treated similarly to the previous consideration. In the 
original Bernoulli game, the lotteries $L_n$ are characterized by $U_n=n$ increasing 
to infinity with rising $n$. However the stochastically preferred is a lottery 
$L_{n_{opt}}$ with a finite $n=n_{opt}$ given by $1/|\beta|$. For repeated games, the 
compound lotteries $L_N$ of $N$ repeated games correspond to the values $U_N=1+\log_2N$ 
increasing with increasing $N$. However, the optimal lottery, with the optimal number 
of runs $N_{opt}$, is described by a finite $N=N_{opt}$ given by (\ref{45}).

\section{Inverse St. Petersburg Paradox}

There exists a situation that in some sense is inverse to the St. Petersburg paradox, 
when one plays a game resulting in losses, whose expected values or expected utilities 
are negative and tend to $-\infty$. Nevertheless, some people, of course not all but 
maybe a small fraction, still continue playing. A typical example concerns a roulette. 
Let us keep in mind the often used roulette whose wheel has $38$ numbers, which consist 
of the numbers one through $36$, the number zero, and the number double zero. The zero 
and double zero are losing numbers for the players. Thus the odds are in favor of the 
casino when betting on roulette. Half of the numbers between $1$ and $36$ are red, and 
the other half are black. It is possible to place a bet just on a number, on a number 
being red or black, odd or even, first $18$ or last $18$, and other combinations.

There is a well-known so-called martingale strategy, used by gamblers which appears 
to guarantee a win. Every time a player loses, he/she doubles his/her bet on the next 
spin of the wheel and continuing to do so on each successive loss, when eventually
he/she does win. If so, he/she will win back the amount he/she lost in addition to a 
profit equal to the original stake. However, this strategy leads to consecutive losses 
rising to infinity. Then the question is: How would it be possible to explain this 
paradox, when some people keep gambling and losing, instead of stopping?

As an illustration, we can take the example considered by Aloysius (2003), when the 
player bets on colour or on first/last $18$ numbers. The process is assumed to terminate 
when the gambler wins. The expected value at the $n$-th stage, after the $n$-th spin 
of the wheel, is 
\be
\label{46}
U_n = \sum_{k=0}^{n-1} p ( 1 - p )^k x_0 \; + \; ( 1 - p )^n x_n \;  ,
\ee 
where $x_0$ is the initial bid, $x_n$ is the loss accumulated after the $n$-th step, 
and
\be
\label{47}
 x_n = - \sum_{k=0}^{n-1} \; 2^k \; x_0 \; , \qquad p = \frac{18}{38} \;  .
\ee
The values of $U_n$ are negative and decreasing, so that
\be
\label{48}
 \lim_{n\ra\infty} \; U_n = - \infty \;  .
\ee
Aloysius (2003) shows that it is possible to find such utility functions $u(x)$ that 
the expected utility
\be
\label{49}
\widetilde U_n = 
\sum_{k=0}^{n-1} \; p ( 1 - p )^k \; u( x_0 ) \; + \; ( 1 - p )^n \; u( x_n ) 
\ee
enjoys the property
\be
\label{50}
  \lim_{n\ra\infty} \; \widetilde U_n ~ \geq ~ u( 0 ) \; .
\ee
Then for the people characterized by such special utility functions, it looks 
reasonable to continue playing.

This explanation is similar to that used for making the expected utility finite by 
introducing concave utility functions in Sec. 3.1 when considering the St. Petersburg 
paradox. However, as has been explained in that section, it is easy to change the 
payoffs so that the problem reduces back to the initial formulation with the infinite 
expected utility. In the same way, here it is possible to change the payoffs $x$ to 
$y(x)$ such that $u(y(x))=x$. Then expression (\ref{49}) returns back to (\ref{46}) 
and the problem reappears. As before, it is interesting to pose a more general question: 
Why some people continue gambling despite the increasing losses tending to infinity?

The general answer is as follows. The case of gambling is different from the choice 
between the lotteries in the St. Petersburg paradox. In the latter, when hesitating 
what ticket to buy, one considers a set of possible lotteries $L_n$, choosing the 
most preferable among them. In the martingale roulette game, one never considers all 
possibilities at once, imagining numerous forthcoming losses, but at each step the 
gambler compares just two alternatives: to stop or keep playing. This choice repeats 
at each next step. The probability of stopping or playing is given by the probability 
$p(L_n)$, as discussed in Sec. 4. Suppose one chooses between the alternatives $L_n$ 
and $L_{n+1}$. Despite that the corresponding expected utilities are negative and 
decreasing, so that $U_{n+1}<U_n$, the probabilities $p_1(L_n)$ (stop gambling) and 
$p_2(L_{n+1})$ (continue gambling) are finite, which means that at each step there 
are some people deciding to stop, as well as there is a fraction of subjects deciding 
to keep gambling.         

Let us take the expected value (\ref{46}) that equals the expected utility under a 
linear utility function. This expected value, by summing the terms, can be rewritten 
as  
\be
\label{51}
 U_n = \left[\; 1 - 2^n ( 1 - p )^n \; \right] \; x_0 \;  .
\ee
For $p = 18/38$, this takes the form
\be
\label{52}
 U_n = \left[\; 1 - \left( \frac{20}{19} \right)^n \; \right] \; x_0 \;  .
\ee
After the first spin of the roulette wheel, there is the choice either to stop playing, 
having the expected value $U_1 = -0.0526 x_0$ or to continue gambling getting the value 
$U_2 = -0.108 x_0$. The probabilities of alternatives with negative expected utilities 
are given by expression (\ref{20}). For simplicity, let us accept the case of neutral 
beliefs, when $\bt=0$. For finite $\beta$, the situation is qualitatively similar. Then 
we find $p_1(L_1)=0.671$ and $p_2(L_2)=0.328$. This implies that, despite the loss, 
about $33\%$ of gamblers continue gambling.    

After the second spin of the wheel, one needs to compare the alternatives $L_2$, with 
$U_2=-0.108 x_0$, and $L_3$, with $U_3=-0.166 x_0$, which gives $p_1(L_2)=0.606$ and 
$p_2(L_3)=0.394$. At the next stage, one compares the alternatives $L_3$ and $L_4$,
with the values $U_3=-1.166x_0$ and $U_4=-0.228x_0$, which yields $p_1(L_3)=0.579$ and 
$p_2(L_4)=0.421$. Then, comparing the alternatives $L_4$ and $L_5$, with $U_4=-0.228x_0$ 
and $U_5=-0.292x_0$, one gets  $p_1(L_4)=0.562$ and $p_2(L_5)=0.438$. For large $n$, 
such that $n\gg 40$, one has
\be
\label{53}
  U_n \simeq  - \left( \frac{20}{19} \right)^n  \; x_0 \;  .
\ee
Then the probabilities of stopping or continuing the roulette gamble tend to
$$
 p_1(L_n) = 0.513 \; , \qquad p_2(L_{n+1}) = 0.487 \qquad ( n \gg 40 ) \; .
$$
This means that about $49\%$ of people, who have reached the $n$-th step, decide to 
keep gambling despite permanent losses.

\section{Conclusion}

In order to resolve the Bernoulli St. Petersburg paradox, there is no need to change
the definitions of either utility functions or expected utilities. The fact that the
expected utility of the Bernoulli game, in the limit, can be infinite is unavoidable.
The Bernoulli game is mathematically correct and defines a bone-fide lottery.

The paradox arises because the lottery with a maximal expected utility corresponds to
infinite price, although the majority of real human beings prefer a lottery that does
not cost more than one or two tosses. This means that people does not consider as
optimal a lottery with a maximal expected utility. That is, the paradox is in the
incompatibility of the definition of an optimal lottery as that for which the expected
utility is maximal and the behavior of real humans.

However, in the frame of the probabilistic approach, defining a stochastically optimal
lottery as that one corresponding to the maximal probability distribution, there is no
any paradox, just because a lottery with infinite expected utility is not stochastically
optimal. On the contrary, a stochastically optimal is a lottery with a rather low
utility corresponding, for real humans, to a couple of tosses.

\section*{Acknowledgment}

The author is grateful for discussions to D. Sornette and E.P. Yukalova.

\vskip 2mm

This research did not receive any specific grant from funding agencies in the public,
commercial, or not-for-profit sectors.

\newpage

{\Large {\bf References}}

\vskip 5mm

{\parindent=0pt

\vskip 2mm
Aloysius, J.A. (2003):
``Rational escalation of costs by playing a sequence of unfavorable gambles: the martingale",
{\it Journal of Economic Behavior and Organization}, 51 (2003) 111--129.

\vskip 2mm
Arieli, A., A. Sterkin, A. Grinvald, and A. Aertsen (1986):
``Dynamics of ongoing activity: explanation of the large variability in evoked cortical 
responses", 
{\it Science}, 273, 1868--1871.

\vskip 2mm
Bernoulli, D. (1738):
``Exposition of a new theory on the measurement of risk",
{\it Proceedings of the Imperial Academy of Sciences of St. Petersburg}, 5, 175--192.
Reprinted in (1954): {\it Econometrica}, 22, 23--36.

\vskip 2mm
Blavatskyy, P. (2005):
``Back to the St. Petersburg paradox?"
{\it Management Science}, 51, 677--678. 

\vskip 2mm
Blavatskyy, P. (2006):
``Violations of betweenness or random errors?",
{\it Economics Letters}, 91, 34--38.

\vskip 2mm
Bottom, W.P., R.N. Bontempo, and  D.R. Holtgrave (1989):
``Experts, novices, and the St. Petersburg Paradox: Is one solution enough?",
{\it Journal of Behavioral Decision Making}, 2, 113--121.

\vskip 2mm
Burnham, K.P., and D.R. Anderson (2002):
{\it Model Selection and Multi-Model Inference}. Berlin: Springer. 

\vskip 2mm
Carbone, E. (1997):
``Investigation of stochastic preference theory using experimental data",
{\it Economics Letters}, 57, 305--311.

\vskip 2mm
Clark, M. (2002):
``The St. Petersburg paradox", in {\it Paradoxes from A to Z}.
London: Routledge, p. 174--177.

\vskip 2mm
Cover, T.M., and J.A. Thomas (2006):
{\it Elements of Information Theory}. New York: Wiley.

\vskip 2mm
Cox, J.C., V. Sadiraj, and B. Vogt (2009):
``On the empirical relevance of St. Petersburg lotteries",
{\it Economics Bulletin}, 29, 214--220.

\vskip 2mm
Cox, J.C., E.B. Kroll, M. Lichters, V. Sadiraj, and B. Vogt (2018):
``The St. Petersburg paradox despite risk-seeking preferences: an experimental study",
{\it Business Research}, 12, 27--44.

\vskip 2mm
Cramer, G. (1728):
``Letter to Nicolas Bernoulli", London, 21 May 1728. \\
(\url{http://www.cs.xu.edu/math/Sources/Montmort/stpetersburg.pdf#search=%22Nicolas%20
Bernoulli%22}).

\vskip 2mm
Dmitruk, A.V., and N.V. Kuzkina (2005):
``Existence theorem for optimal control problems on an infinite time interval",
{\it Mathematical Notes}, 78, 466--480.

\vskip 2mm
Fechner, G. (1860):
{\it Elements of Psychophysics}. New York: Holt, Rinehart and Winston.

\vskip 2mm
French, S., and D.R. Insua (2000):
{\it Statistical Decision Theory}. London: Arnold. 

\vskip 2mm
Gelfand, I.M., and S.V. Fomin (1963):
{\it Calculus of Variations}. London: Prentice Hall.

\vskip 2mm
Gerken, J. (1991):
``Generalized logit model",
{\it Transportation Research B}, 25, 75--88.

\vskip 2mm
Ghahramani, S. (2000):
{\it Fundamentals of Probability}. New Jersey: Prentice Hall. 

\vskip 2mm
Glimcher, P.W. (2005):
``Indeterminacy in brain and behavior", 
{\it Annual Review of Psychology}, 56, 25--56.

\vskip 2mm
Glimcher, P.W., and A.A. Tymula (2018):
``Expected subjective value theory (ESVT): A representation of decision under risk and 
certainty",
University of Sydney Economics Working Paper Series.  

\vskip 2mm
Gold, J.I., and M.N. Shadlen (2001):
``Neural computations that underlie decisions about sensory stimuli", 
{\it Trends in Cognitive Science}, 5, 10--16.

\vskip 2mm
Gul, F., P. Natenzon, and W. Pesendorfer (2014):
``Random choice as behavioral optimization",
{\it Econometrica}, 82, 1873--1912.

\vskip 2mm
Gustason, W. (1994):
{\it Reasoning from Evidence}. New York: Macmillan College.

\vskip 2mm
Hardin, R. (1982):
{\it Collective Action}. Baltimore: John Hopkins University.

\vskip 2mm
Harless, D., and C. Camerer (1994):
``The predictive utility of generalized expected utility theories",
{\it Econometrica}, 62, 1251--1289.

\vskip 2mm
Hayden, B.Y., and M.L. Platt (2009):
``The mean, the median, and the St. Petersburg paradox”,
{\it Judgment and Decision Making}, 4, 256--272.

\vskip 2mm
Hey, J.D., and C. Orme (1994):
``Investigating generalisations of expected utility theory using experimental data",
{\it Econometrica}, 62, 1291--1326.

\vskip 2mm
Jeffrey, R.C. (1983):
{\it The Logic of Decision}. Chicago: University of Chicago.

\vskip 2mm
Kahneman, D., and A. Tversky (1979):
``Prospect theory: an analysis of decision under risk",
{\it Econometrica}, 47, 263--291.

\vskip 2mm
Klyve, D., and A. Lauren (2011):
``An empirical approach to the St. Petersburg paradox",
{\it College Mathematics Journal}, 42, 260--264.

\vskip 2mm
Kullback, S., and R.A. Leibler (1951):
``On information and sufficiency",
{\it Annals of Mathematical Statistics}, 22, 79--86.

\vskip 2mm
Kullback, S. (1959):
{\it Information Theory and Statistics}. New York: Wiley.

\vskip 2mm
Kurtz-David, V., D. Persitz, R. Webb, and D.J. Levy (2019):
``The neural computation of inconsistent choice behavior",
{\it Nature Communications}, 10, 1583.

\vskip 2mm
Loomes, G., and R. Sugden (1995):
``Incorporating a stochastic element into decision theories",
{\it European Economic Review}, 39, 641--648.

\vskip 2mm
Loomes, G., P. Moffatt, and R. Sugden (2002):
``A microeconomic test of alternative stochastic theories of risky choice",
{\it Journal of Risk and Uncertainty}, 24, 103--130.

\vskip 2mm
Luce, R.D. (1958):
``A probabilistic theory of utility",
{\it Econometrica}, 26, 193--224.

\vskip 2mm
Luce, R.D. (1959):
{\it Individual Choice Behavior: A Theoretical Analysis}. New York: Wiley.

\vskip 2mm
Luce, R.D. (1977):
``The choice axiom after twenty years",
{\it Journal of Mathematical Psychology}, 15, 215--233.
 
\vskip 2mm
MacKay, D.J.C. (2003):
{\it Information Theory, Inference, and Learning Algorithms}. 
Cambridge: Cambridge University Press.

\vskip 2mm
Martin, R.M. (2008):
``The St. Petersburg paradox", in {\it Stanford Encyclopedia of Philosophy}.
Stanford University: Stanford.

\vskip 2mm
Menger, K. (1967):
``The role of uncertainty in economics", in
{\it Essays in Mathematical Economics in Honor of Oscar Morgenstern},
edited by M. Shubik. Princeton: Princeton University.

\vskip 2mm
Moore, A.W. (1990):
{\it The Infinite}. London: Routledge.

\vskip 2mm
Nathan, A. (1984):
``False expectations", {\it Philosophy of Science}, 51, 128--136.

\vskip 2mm
Neugebauer, T. (2010):
``Moral impossibility in the St Petersburg paradox: A literature survey and experimental
evidence",
{\it Luxembourg School of Finance Working Paper}, 1--43.

\vskip 2mm
Nielsen, R. (2005):
{\it Statistical Methods in Molecular Evolution}. Berlin: Springer.
 
\vskip 2mm
Nobandegani, A.S., and T.R. Shultz (2020):
``The St. Petersburg paradox: A fresh algorithmic perspective",
{\it Association for the Advancement of Artificial Intelligence Working Paper}, 1--4.

\vskip 2mm
Pfiffermann, M. (2011):
``Solving the St. Petersburg paradox in cumulative prospect theory: the right amount
of probability weighting",
{\it Theory and Decision}, 71, 325--341.

\vskip 2mm
Pinski, F.J., G. Simpson, A.M. Stuart, and H. Weber (2015):
``Kullback-Leibler approximation for probability measures on infinite dimensional spaces",
{\t SIAM Journal of Mathematical Analysis}, 47, 4091--4122.

\vskip 2mm
Poincar\'e, H. (1902):
{\it Science and Hypothesis}. London: Walter Scott.

\vskip 2mm
Rayo, L., and G.S. Becker (2007):
``Evolutionary efficiency and happiness",
{\it Journal of Political Economy}, 2007, 115, 302--337.

\vskip 2mm
Rieger, M.O., and M. Wang (2006):
``Cumulative prospect theory and the St. Petersburg paradox",
{\it Economic Theory}, 28, 665--679.

\vskip 2mm
Rivero, J.C., D.R. Holtgrave, R.N. Bontempo, and W.P. Bottom (1990):
``The St. Petersburg Paradox: Data, at last".
{\it Commentary}, 8, 46--51.

\vskip 2mm
Resnik, M.D. (1987):
{\it Choices: An Introduction to Decision Theory}.
Minneapolis: University of Minnesota.

\vskip 2mm
Rubinstein, R.Y., and D.P. Kroese (2004):
{\it The Cross-Entropy Method}. Berlin: Springer.

\vskip 2mm
Rucker, R. (1995):
{\it Infinity and the Mind: The Science and Philosophy of the Infinite}.
Princeton: Princeton University.

\vskip 2mm
Samuelson, P.A. (1960):
``The St. Petersburg paradox as a divergent double limit",
{\it International Economic Review}, 1, 31--37.

\vskip 2mm
Samuelson, P.A. (1977):
``St. Petersburg paradoxes: defanged, dissected, and historically described",
{\it Journal of Economic Literature}, 15, 24--55.

\vskip 2mm
Savage, L. (1954):
{\it The Foundations of Statistics}. New York: Wiley.

\vskip 2mm
Schumacher, J.F., S.K. Thompson, and C.A. Olman (2011):
``Contrast response functions for single Gabor patches: ROI-based analysis over-represents 
low-contrast patches for GE BOLD", 
{\it Frontiers in System Neuroscience} 5, 1--10. 

\vskip 2mm
Seidl, C. (2013):
``The St. Petersburg Paradox at 300",
{\it Journal of Risk and Uncertainty}, 46, 247--264.

\vskip 2mm
Shadlen, M.N. and D. Shohamy (2016): 
``Perspective decision making and sequential sampling from memory", 
{\it Neuron}, 90, 927--939.

\vskip 2mm
Shore, J., and R. Johnson (1980):
``Axiomatic derivation of the principle of maximum entropy and the principle of
minimum cross-entropy",
{\it IEEE Transactions on Information Theory}, 26, 26--37.

\vskip 2mm
Slovic, P., C.K. Mertz, D.M. Markowitz, A. Quist, and D. V\"{a}stfj\"{a}ll (2020): 
``Virtuous violence from the war room to death row",
{\it Proceedings of the National Academy of Sciences of USA} 117, 20474--20482. 

\vskip 2mm
Small, K. (1987):
``A discrete choice model for ordered alternatives",
{\it Econometrica}, 55, 409--424. 

\vskip 2mm
Tversky, A., and M. Bar-Hillel (1983):
``Risk: The long and the short",
{\it Journal of Experimental Psychology: Learning, Memory, and Cognition},
9, 713--717 

\vskip 2mm
Vivian, R.W. (2004):
``Simulating Daniel Bernoulli’s St. Petersburg game: Theoretical and empirical
consistency”.
{\it Simulation and Gaming}, 35, 499--504.

\vskip 2mm
Vivian, R.W. (2013):
``Ending the myth of the St. Petersburg paradox",
{\it South African Journal of Economic and Management Sciences}, 16, 347--364.

\vskip 2mm
von Neuman, J., and O. Morgenstern (1944):
{\it Theory of Games and Economic Behavior}. Princeton: Princeton University.

\vskip 2mm
Webb, R. (2019):
``The neural dynamics of stochastic choice", 
{\it Management Science}, 65, 230--255.

\vskip 2mm
Weirich, P. (1984):
``The St. Petersburg gamble and risk",
{\it Theory and Decision}, 17, 193--202.

\vskip 2mm
Wen, C.H., and F.S. Koppelman (2001):
``The generalized nested logit model",
{\it Transportation Research B}, 35, 627--641.

\vskip 2mm
Werner, G., and V.B. Mountcastle (1963):
``The variability of central neural activity in a sensory system and its implications 
for the central reflection of sensory events", 
{\it Journal of Neurophysiology}, 26, 958--977.

\vskip 2mm
Williams, R. (2016):
``Understanding and interpreting generalized ordered logit models",
{\it Journal of Mathematical Sociology}, 40, 7--20. 

\vskip 2mm
Woodford, M. (2020):
``Modeling imprecision in perception, valuation, and choice",
{\it Annual Review of Economics}, 12, 579--601. 

\vskip 2mm
Yukalov, V.I., and D. Sornette (2014):
``Role of information in decision making of social agents",
{\it International Journal of Information Technology and Decision Making}, 14, 1129--1166.

\vskip 2mm
Yukalov, V.I., and D. Sornette (2018):
``Quantitative predictions in quantum decision theory",
{\it IEEE Transactions on Systems, Man and Cybernetics Systems}, 48, 2168--2216.

}

\end{document}